\documentclass[12pt,reqno]{article}

\usepackage{amsfonts}
\usepackage{amssymb}
\usepackage{amsbsy}
\usepackage{amsmath}

\title{Almost giant  clusters  for percolation
  \\
on large trees with logarithmic heights} 
\author{Jean Bertoin\thanks{Institut f\"ur Mathematik, 
Universit\"at Z\"urich, 
Winterthurerstrasse 190, 
CH-8057 Z\"urich, Switzerland. \hfill \eject
Email: jean.bertoin@math.uzh.ch} }

\date{}

\usepackage{graphicx}
\usepackage{amsfonts,amsmath,amssymb,bbm}
\usepackage{setspace} 
\def\proof{\noindent{\bf Proof:}\hskip10pt}        
\def\QED{\hfill $\Box$}

\font\tenmath=msbm10 scaled 1200
\font\sevenmath=msbm7 scaled 1200
\font\Fivemath=msbm5 scaled 1200
\newfam\mathfam \textfont\mathfam=\tenmath
\scriptfont\mathfam=\sevenmath \scriptscriptfont\mathfam=\Fivemath

\def \\ { \cr }
\def\R{\mathbb{R}}
\def \1{1 \mkern -6mu 1} 
\def\N{\mathbb{N}}
\def\E{\mathbb{E}}
\def\P{\mathbb{P}}

\def\R{\mathbb{R}}

\def \e{{\rm e}}
\def \x{{\bf x}}

\def \d{{\rm d}}

\newtheorem{theorem}{Theorem}

\newtheorem{proposition}{Proposition}

\newtheorem{corollary}{Corollary}
\begin{document}

\maketitle

\begin{abstract} 
This text is based on a lecture for the Sheffield Probability Day; its main purpose is to survey some recent asymptotic results \cite{Be, BU} about  Bernoulli bond percolation on certain large random trees with logarithmic height. We also provide a general criterion for the existence of giant percolation clusters in large trees, which answers a question raised by David Croydon.
\end{abstract}

{\bf Key words:} Random tree; percolation; giant component.

{\bf Subject Classes:} 60K35; 05C05.
  
\begin{section}{Introduction}

 It is well-known that percolation is considerably simpler to study on a tree than on a general graph, thanks to the property of
 uniqueness of the path connecting two vertices. We refer in particular to Chapter 5 of \cite{LP} and references therein for a number of  important and useful results  for infinite trees, such as criteria for the existence or absence of infinite percolation clusters. 
Here, we shall be interested in a somewhat  different type of questions. Specifically we consider a tree of large but finite size, perform  a Bernoulli bond percolation with a parameter that depends on the size  of that tree, and our purpose is to investigate the asymptotic behavior of  the sizes of the largest clusters  for appropriate regimes when the size of the tree goes to infinity. 

Our motivation comes from a celebrated result of Erd\H{o}s and R\'enyi on the random graph model, which can be phrased informally as follows. With  high probability when $n\gg 1$, Bernoulli bond percolation on the complete graph with $n$ vertices and with  parameter $p(n)\sim c/n$ for some fixed $c>1$,
produces a single giant cluster of size close to $\theta(c)n$, where $\theta(c)\in(0,1)$ is some  known constant, while the second, third, etc. largest clusters are almost microscopic, and more precisely have size  only of order $\ln n$. 

In the first part of this text, we provide a simple characterization of tree families and percolation regimes which yield giant clusters, 
answering a question raised by David Croydon. In the second part, we  review briefly the main results of \cite{Be, BU} concerning two natural families of random trees with logarithmic heights, namely recursive trees and scale-free trees. 
In those works, we show that the next largest clusters are almost giant, in the sense that their sizes are of order $n/\ln n$, and obtain precise limit theorems in terms of certain Poisson random measures. 
A common feature in the analysis of percolation for these models is that, even though one addresses a static problem, it is useful to consider dynamical versions in which edges are removed, respectively vertices are inserted, one after the other in certain order as time passes. 

\end{section}

\begin{section}{Giant clusters}
We first introduce notations and hypotheses which will have an important role in this section. 
For a given integer $n$, we consider a set  of  $n+1$ vertices, say  ${\mathcal V}_n=\{0,1,\ldots, n\}$, and a tree structure $T_n$ on ${\mathcal V}_n$. So $T_n$ has  $n$ edges, and we should  think of $0$ as the root of $T_n$.
We perform a Bernoulli bond percolation on $T_n$ with parameter $p(n)$, so that each edge of $T_n$ is kept with probability $p(n)$ and removed with probability $1-p(n)$, independently of the other edges. The resulting connected components are then referred to as clusters. 

We write $C_{p(n)}^0$ for the size of the cluster that contains the root; plainly $C^0_{p(n)}\leq n+1$. 
We say that $C_{p(n)}^0$ is {\it giant} if  $n^{-1}C^0_{p(n)}$  converges in law to some random variable $G\not \equiv 0$, which should be thought of as  the asymptotic proportion of vertices pertaining  to the root cluster.
David Croydon raised the question of finding a simple criterion for $C_{p(n)}^0$ to be giant, 
depending of course on the nature  of $T_n$ and regimes of the percolation parameter $p(n)$; 
this motivates the following.

For each fixed $n\in\N$, we denote by $V_1, V_2, \ldots$ a sequence of i.i.d. vertices in ${\mathcal V}_n$ with the uniform distribution. Next, for every $k\in\N$, we write
$L_{k,n}$ for the length of the tree $T_n$ reduced to $V_1, \ldots, V_k$ and the root $0$, i.e. the minimal number of edges of $T_n$ which are needed to connect $0$ and $V_1, \ldots, V_k$. In particular, $L_{1,n}$ should be thought of as  the height of a typical vertex in $T_n$. 
Let $\ell: \N\to \R_+$ be some function with $\lim_{n\to \infty}\ell(n)=\infty$. We introduce the hypothesis
\begin{equation} \label{Hk}
 \frac{1}{\ell(n)}\, L_{k,n}Ê\Rightarrow {L}_k\,,  \tag{$H_k$}
 \end{equation}
where  $\Rightarrow$ means  weak convergence and ${L}_k$ is some random variable with values in $\R_+$. 
We stress that \eqref{Hk} can be assumed to hold for different values of $k$ and then only convergence in the sense of one-dimensional distributions is involved.

In several examples, the function $\ell$ is a logarithm  and ${L}_k\equiv a k$ with $a$ a positive constant.  For instance, this happens for some important families of random trees, such as recursive trees, binary search trees, etc.; see \cite{Drmota}.  Aldous \cite{Al} considered a different class of examples, including the case  when $T_n$ is a Cayley tree of size $n+1$ (i.e. a tree picked uniformly at random amongst the $(n+1)^{n-1}$ trees on ${\mathcal V}_n$),  for which it is known that  \eqref{Hk} holds with $\ell(n)=\sqrt n$ and ${L}_k$ a chi-variable with $2k$ degrees of freedom.

We now state the central result of this section.

\begin{theorem} \label{T1} For an arbitrary $c\geq 0$, consider 
 the regime
\begin{equation}\label{E1}
p(n) = 1-\frac{c}{\ell(n)} + o(1/\ell(n)).
\end{equation}

\noindent {\rm (i)} If \eqref{Hk} holds for every $k\in\N$,
then we have in the regime \eqref{E1}
\begin{equation} \label{E2}
n^{-1}C_{p(n)}^0\Rightarrow {G(c)}\,,
\end{equation}
where ${G(c)}\not \equiv0$ is a random  variable whose law  is determined by its entire moments:
\begin{equation} \label{E4}
\E( {G(c)}^k)= \E(\e^{-c {L}_k})\,,\qquad k\in\N\,.
\end{equation} 
In particular $\lim_{c\to 0+}{G(c)}=1$ in probability. 

\noindent {\rm (ii)} Conversely,  suppose that  for every $c>0$,  \eqref{E2} holds in the regime \eqref{E1}
for some random variable ${G(c)}$ with values in $ [0,1]$. Suppose further that $\lim_{c\to 0+}{G(c)}=1$ in probability. Then  \eqref{Hk} is fulfilled  for every $k\geq 1$,
with ${L}_k$ a nonnegative random variable whose Laplace transform is given by \eqref{E4}.

\end{theorem}

\proof  The proof relies on the observation that for each $k\geq 1$, there is the identity
\begin{equation}\label{E3}
\E\left( \left((n+1)^{-1}ÊC_{p(n)}^0\right)^k\right) = \E\left( p(n)^{L_{k,n}}\right).
\end{equation}
Indeed, recall that $V_1, \ldots, V_k$ are $k$ i.i.d. uniformly distributed vertices, which  are independent of the percolation process.  
This enables us to interpret the left-hand side of \eqref{E3} as the probability that $V_1, \ldots, V_k$ belong to the percolation cluster containing the root.
On the other hand, considering the tree reduced to $V_1, \ldots, V_k$ and the root shows that this same probability can also be expressed in terms of the length $L_{k,n}$ of this reduced tree, as
the right-hand side of \eqref{E3}. 

The assumption \eqref{Hk} entails that in the regime \eqref{E1}, 
$$\lim_{n\to\infty} \E\left( p(n)^{L_{k,n}}\right) = \lim_{n\to \infty}\E\left( \exp\left( - \frac{c}{\ell(n)}  L_{k,n}\right)\right)  =  \E(\e^{-c {L}_k})\,,$$
 and then we deduce from \eqref{E3} that
$$
\lim_{n\to \infty} \E\left( \left((n+1)^{-1}ÊC_{p(n)}^0\right)^k\right) = \E(\e^{-c {L}_k})\,.
$$
Thus,  if \eqref{Hk} holds for every $k\in\N$,  then 
$(n+1)^{-1}ÊC_{p(n)}^0$ converges in law to some variable ${G(c)}$ with values in $[0,1]$. More precisely, the law of ${G(c)}$ is determined by its entire moments $\E( {G(c)}^k)= \E(\e^{-c {L}_k})>0$; this proves \eqref{E2}.

Conversely, as the size of a cluster cannot exceed $n+1$, \eqref{E2} implies that in the regime \eqref{E1}, we have for every integer $k\geq 1$
$$\lim_{n\to \infty} \E\left( \left((n+1)^{-1}ÊC_{p(n)}^0\right)^k\right) = \E( {G(c)}^k)\,.$$
From \eqref{E3}, we rewrite this as 
$$\lim_{n\to \infty} \E\left( p(n)^{L_{k,n}}\right) = \E( {G(c)}^k)\,.$$
Plugging the expression \eqref{E1} for the parameter $p(n)$, we easily derive  that 
$$\lim_{n\to \infty}\E\left( \exp\left( - \frac{c}{\ell(n)}  L_{k,n}\right)\right)  = \E( {G(c)}^k)\,.$$
Recall the assumption that $\lim_{c\to 0+}{G(c)}=1$ in probability, in particular $\lim_{c\to 0+}\E( {G(c)}^k)=1$. 
We conclude from Theorem XIII.1.2 in Feller \cite{Feller} on page 431, that for each $k\geq 1$, the function $c\mapsto \E( {G(c)}^k)$ is the Laplace transform of a random variable ${L}_k\geq 0$, and that 
\eqref{Hk} holds.  \QED

We next point at an interesting consequence of Theorem \ref{T1} to the characterization of the cases for which the proportion of vertices in the root cluster converges in probability to a constant.
We consider the situation where the variables ${L}_k$ appearing in Hypotheses \eqref{Hk} are of the form 
\begin{equation}\label{H'k}
{L}_k=\xi_1+\cdots+ \xi_k \tag{$H'_k$}
\end{equation}
where $\xi_1, \ldots$ is a sequence of i.i.d. variables in $\R_+$.

\begin{corollary} \label{C1} 
\noindent {\rm (i)} 
Suppose that \eqref{Hk} and \eqref{H'k} hold for $k=1,2$.  Then in the regime \eqref{E1}, 
we have 
\begin{equation} \label{E2'}
\lim_{n\to \infty} n^{-1}ÊC_{p(n)}^0 = \theta(c) \qquad \hbox{in probability}\,,
\end{equation}
where $\theta(c) = \E(\e^{-c \xi_1})>0$. Further \eqref{Hk} and \eqref{H'k} hold for every $k\geq 1$. 

\noindent {\rm (ii)} Conversely,  suppose that  for every $c>0$,
 \eqref{E2'} holds in the regime \eqref{E1} for some function $\theta: [0,\infty)\to [0,1]$ such that $\lim_{c\to 0+}\theta(c)=1$. Then $\theta$ is the Laplace transform of a nonnegative random variable $\xi$,
and \eqref{Hk} and \eqref{H'k} are  fulfilled  for every $k\geq 1$ with  $\xi_1, \ldots$ a sequence of i.i.d. copies of $\xi$.

\end{corollary}

\proof  When \eqref{H'k} holds, we have $\E(\exp(-c{L}_k))=\theta(c)^k$, with $\theta(c)=\E(\e^{-c\xi_1})$. 
We now see from the proof of Theorem \ref{T1} that 
Hypotheses \eqref{Hk} and \eqref{H'k} entail that in the regime \eqref{E1}, we have
$$\lim_{n\to \infty} \E\left( \left(n^{-1}ÊC_{p(n)}^0\right)^k\right) = \theta(c)^k\,.$$
In particular,  if \eqref{Hk} and \eqref{H'k} hold for $k=1,2$,  then 
$$\lim_{n\to \infty} \E\left( \left( n^{-1}ÊC_{p(n)}^0- \theta(c)\right)^2\right) =0\,,$$
which proves \eqref{E2'}.

Conversely, if \eqref{E2'} holds, then we can apply Theorem \ref{T1}(ii) with ${G(c)}\equiv \theta(c)$. In particular we know that \eqref{Hk} holds for all $k\in\N$.
Further, we get that $\E(\e^{-c{L}_k})=\theta(c)^k$, which in turn shows  that \eqref{H'k} is fulfilled. \QED
 
We now conclude this section by pointing at a simple criterion which ensures that the cluster containing the root is the unique giant component. 

\begin{proposition}\label{P1}
In the preceding notation, assume that there is the joint weak convergence
$$ \frac{1}{\ell(n)}( L_{1,n}, L_{2,n}) Ê\Rightarrow (L_1, L_2)\,,$$
where $(L_1, L_2)$ is a pair of random variables such that $L_2-L_1$ has the same law as $L_1$.
Then for every $c>0$,  in the regime \eqref{E1}, we have
$$\lim_{n\to \infty} n^{-1}ÊC_{p(n)}^1 = 0 \qquad \hbox{in probability,}$$
where $C^1_{p(n)}$ denotes  the size of the largest percolation cluster which does not contain the root $0$. 
\end{proposition}

\proof Recall that $V_1$ and $V_2$ denote two independent uniformly distributed random vertices. 
Plainly,   the probability $ \varrho(n)$ that $V_1$ and $V_2$ both belong to the same percolation cluster and are disconnected from $0$ can be bounded from below by
$(n+1)^{-2}\E(|C^1_{p(n)}|^2)$

On the other hand, $ \varrho(n)$ is bounded from above by the probability that at least one edge of the branch from the root $0$ to the  branch-point $V_1\wedge V_2$ of $V_1$ and $V_2$ has been removed, viz.
\begin{equation}\label{E5}
(n+1)^{-2}\E(|C^1_{p(n)}|^2) \leq  \varrho(n) \leq 1-\E\left( p(n)^{ \d_n(0,V_1\wedge V_2)} \right), 
\end{equation}
where $\d_n$ denotes  the graph distance in $T_n$. 
 
 Next, write $$L_{2,n}=\d_n(0,V_1)+\d_n(0,V_2)-\d_n(0,V_1\wedge V_2)\,.$$
Since 
 $L_{1,n}=\d_n(0,V_1)$ has the same law as $\d_n(0,V_2)$, it follows  from our assumption that the 
 sequences $\ell(n)^{-1} \d_n(0,V_2)$ and $\ell(n)^{-1}\left( \d_n(0,V_2)-\d_n(0,V_1\wedge V_2)\right)$ 
  converge weakly to the same distribution. This readily implies that
  $$
\d_n(0,V_1\wedge V_2)= o(\ell(n))\qquad \hbox{in probability,}$$
and we conclude that the right-hand side in \eqref{E5}  tends to $0$ as $n\to \infty$. \QED

\end{section}

\begin{section}{Almost giant clusters}
In this section, we turn our attention to the percolation clusters which do not contain the root. 
We  write
$$C^1_{p(n)} \geq C^2_{p(n)} \geq \ldots$$
for the sequence of their sizes, ranked in the decreasing 
order\footnote{Beware that this convenient  notation may be sightly misleading, since $C^0_{p(n)}$ is always the size of the cluster containing the root $0$,
while  for $i\geq 1$, 
$C^i_{p(n)}$ is in general  not the size of the cluster containing the vertex $i$.}.  A natural problem is then  to determine the asymptotic behavior of this sequence. 
We first point out that Hypotheses \eqref{Hk} are insufficient to characterize the latter, by considering three simple examples in which very different behaviors can be observed.

First, imagine that $T_n$ is a star-shaped tree centered at $0$, meaning  that  the root is the unique branching point. Suppose also for simplicity  that there are $\sim n^{1-\alpha}$ branches attached to the root, each of size $\sim n^{\alpha}$, where  $\alpha\in(0,1)$ is some fixed parameter. Then one readily checks that \eqref{Hk} and \eqref{H'k}  hold for every $k\geq 1$ with $\ell(n)=n^{\alpha}$ and ${L}_k=
\xi_1+\cdots + \xi_k$ where the $\xi_i$ are i.i.d.uniformly distributed on $[0,1]$.  It  is further straightforward to see in the regime \eqref{E1}, one has  
$$C^1_{p(n)} \sim  C^2_{p(n)} \sim \ldots \sim C^j_{p(n)} \sim n^{\alpha}$$
 for every fixed $j\in\N$. 

Second, consider the case when $T_n$ is the complete regular $d$-ary tree with height $h$, where $d\geq 2$ is some integer. So there are
$d^j$ vertices at distance $j=0,1,\ldots, h$ from the root and 
$$n=n(h)= d(d^{h}-1)/(d-1).$$
One readily checks that  Hypotheses \eqref{Hk} and \eqref{H'k} hold for every $k\geq 1$ with  $\ell(n)=\ln n$ and $\xi_i\equiv 1/\ln d$. Because the subtree spanned by a vertex at height $j\leq h$ is again a complete regular $d$-ary tree with height $h-j$, we deduce from the preceding section that in the regime \eqref{E1}, the size $C^1_{p(n)}$ of the largest cluster which does not contain the root is close to
$$\e^{-c/\ln d} d^{h-\kappa(h)+1}/(d-1),$$
where $\kappa(h)$ is the smallest height at which an edge has been removed.
Recall that there are  $d(d^j-1)/(d-1)$ edges with height at most $j$,  so the law of $\kappa(h)$ is given by
$$\P(\kappa(h)> j)= p(n)^{d(d^j-1)/(d-1)}\,,\qquad j=1, \ldots, h\,.$$
It follows readily that in the regime \eqref{E1}, the sequence $\left( \kappa(h)- \frac{\ln h}{\ln d}: h\in\N\right)$ is tight. We stress however that this sequence {\it does not  converge in distribution} as $h\to \infty$; more precisely weakly convergent subsequences are obtained provided that the fractional part $\{\frac{\ln h}{\ln d}\}$ converges. 
It follows that the sequence  
$\left( \frac{\ln n}{n} C^1_{p(n)}: n=n(h), h\in\N\right) $ is also tight. It does not converge as $h\to \infty$; however weakly convergent subsequences can be excerpt  provided that $\{\frac{\ln h}{\ln d}\}$ converges.

Third, we recall that in the case of Cayley trees, Pitman \cite{Pi1, PiSF} showed that 
for $1-p(n)\sim c/\sqrt n $ with a fixed $c>0$, the sequence of the sizes of the clusters ranked in decreasing order and renormalized by a factor $1/n$ converges weakly as $n\to \infty$  to a random mass partition which can be described explicitly in terms of a conditioned Poisson measure. It is interesting to observe that in this situation, the number of giant components is unbounded as $n\to\infty$. We stress that the conditions of Proposition \ref{P1} and the hypotheses \eqref{H'k} for $k\geq 2$ fail for Cayley trees.

We shall now study the asymptotic behavior of the sizes of the largest clusters which do not contain the root for two families of random trees with logarithmic heights, i.e. which fulfill \eqref{Hk} with $\ell(n) = \ln n$. In particular, we shall point out that in the regime \eqref{E1}, the largest percolation clusters which do not contain the root  fail to be giant only by a logarithmic factor.

\subsection{Random recursive trees}
A tree on an ordered set of vertices is called {\it recursive} if,  when we agree that the smallest vertex serves as the root, then the sequence of vertices along any branch from the root to a leaf is increasing.  Recursive trees are sometimes also known as increasing trees in the literature; they arise for instance in computer science as data structures, or as simple epidemic models.  

Of course, there is no loss of generality in assuming that  the set of vertices is ${\mathcal V}_n=\{0,1, \ldots, n\}$ (and then $0$ is the root); however other ordered sets may arise naturally in this setting as we shall see. Each recursive tree  on ${\mathcal V}_n$ encodes a permutation of $\{1, \ldots, n\}$ in such a way that the subtrees attached to the root $0$ correspond to the cycles of  the permutation, and this encoding is bijective; see Section 6.1.1 in \cite{Drmota}. In particular, there are $n!$  recursive trees on ${\mathcal V}_n$; we pick one of them uniformly at random and denote it by $T_n$. In other words, $T_n$ can be viewed as a Cayley tree on ${\mathcal V}_n$, subject to the condition that the sequence of vertices along any branch from the root to a leaf is increasing. We stress that, informally,  the conditioning becomes singular as $n\to \infty$. Indeed the geometry of large Cayley trees and large uniform recursive trees are notoriously different; for instance the typical height of the former is of order $\sqrt n$ while that of latter is only of order $\ln n$.

There is an elementary 
algorithm for constructing $T_n$ which is closely related to the so-called Chinese restaurant process (see, e.g. Section 3.1 in Pitman \cite{PiSF}), and hence further points at the connexion with uniform random permutations. For every $i=1, \ldots, n$,  we pick a vertex $U_i$ uniformly at random from $\{0,\ldots, i-1\}$
and independently of the $U_j$ for $j\neq i$. The random tree  induced by the set of edges $\{(i,U_i): i=1, \ldots , n\}$ is then a version of $T_n$.  

Uniform recursive trees fulfill an important {\it splitting property} which is the key to many of their features. Fix an arbitrary $k\in\{1,\ldots, n\}$ and remove the edge between $k$ and its parent $U_k$. This disconnects $T_n$ into two subtrees, say $T$ and $T'$. If we denote by ${\mathcal V}$ (respectively, ${\mathcal V}'$) the sets of vertices of $T$ (respectively, of $T'$), then conditionally on ${\mathcal V}$ and ${\mathcal V}'$, $T$ and $T'$ are two independent uniform recursive trees with respective sets of vertices ${\mathcal V}$ and  ${\mathcal V}'$. This basic property is easy to check,  either directly from the definition, or from the Chinese restaurant construction of $T_n$. 

It is easy to verify that the conditions \eqref{Hk} are fulfilled for all $k\geq 1$ with $\ell(n)=\ln n$ and $L_k\equiv k$; see Section 6.2.5 in \cite{Drmota}. We conclude from the preceding section that in the regime \eqref{E1}, the cluster containing the root $0$ is the unique giant percolation cluster of $T_n$, and more precisely that \eqref{E2'} holds with $\theta(c) = \e^{-c}$. 
The main result of \cite{Be} is that the next largest clusters are almost giant, and more precisely,  one has the following weak limit theorem.

\begin{theorem}\label{T2} Let $T_n$ denote a uniform random recursive tree on $\{0,1, \ldots, n\}$.  For every fixed integer $j$, in the regime \eqref{E1} with $\ell(n)=\ln n$,  
$$\left(\frac{\ln n}{n}C^1_{p(n)}, \ldots, \frac{\ln n}{n}C^j_{p(n)}\right)$$
converges in distribution towards
$$(\x_1, \ldots, \x_j)$$
where $\x_1>\x_2>\ldots$ denotes the sequence of the atoms of a Poisson random measure
on $(0,\infty)$ with intensity $c\e^{-c} x^{-2}\d x\,.$
\end{theorem}

There is an equivalent simple description of the law of the limiting sequence, namely $1/\x_1$, $1/\x_2-1/\x_1, \ldots, 1/\x_j-1/\x_{j-1}$ are i.i.d. exponential variables with parameter $c\e^{-c}$. In particular $\x_j$ has the same distribution as the inverse of 
 a gamma variable with parameter $(j,c\e^{-c})$, and $\lim_{j\to \infty}j\x_j=c \e^{-c} $ in probability. 

The basic idea in \cite{Be} for establishing Theorem \ref{T2} is to relate percolation on a rooted tree $T$ to a random algorithm for the isolation of its root that was 
introduced Meir and Moon. Specifically, following these authors, we can imagine that we remove
an edge in $T$  uniformly at random, disconnecting $T$ into two subtrees. We set aside the subtree which does not contain the root and iterate in an obvious way with the subtree containing the root, until the root is finally isolated. Loosely speaking, we can think of this algorithm as a dynamical version of percolation (i.e. edges are now removed one after the other rather than simultaneously), except that each time an edge is removed, the cluster which does not contain the root is instantaneously frozen, in the sense that only edges  belonging to the cluster that contains the root can be removed.

The upshot of this point of view is that it enables us to use a coupling due to Iksanov and M\"{o}hle \cite{IM}, which, informally, identifies the sequence of the sizes of the frozen subtrees which arise from the isolation of the root algorithm, with the sequence $\eta_1, \eta_2, \ldots, \eta_k$ of i.i.d. variables with distribution
$$\P(\eta = j)= \frac{1}{j(j+1)}\,,\qquad j\in \N\,,$$
at least as long as $\eta_1+\cdots + \eta_k\leq n$. 
In short, this coupling follows from the splitting property of random recursive trees, and the following remarkable fact observed by Meir and Moon \cite{MM}. Imagine that we remove an edge of $T_n$ uniformly at random, and consider the size of the resulting subtree that does not contain the root. Then the latter has the same distribution as $\eta$ conditioned on $\eta\leq n$.  

The coupling of Iksanov and M\"{o}hle enables us to use  Extreme Values Theory and determine the asymptotic behavior of the sizes of these frozen subtrees, jointly with the steps of the algorithm  at which they have appeared. In short, one finds that the largest frozen sub-trees have size of order $n/\ln n$ and a precise limit theorem can be given in terms of the atoms of some Poisson random measure. 
It then remains de-freeze  each of these subtrees by performing an additional Bernoulli percolation with a suitable parameter, to recover the outcome of percolation on $T_n$. 
Roughly speaking, each of these frozen  subtrees can be viewed conditionally on its size as a uniform recursive tree. As a consequence, the additional percolation produces a single relatively giant component of size again of order $n/\ln n$ and further clusters of smaller size $O\left(n/\ln^2n\right)$. 
In particular, the largest percolation clusters of $T_n$ which do not contain the root correspond to simple transformations of the frozen subtrees arising from the algorithm of isolation of the root, and their limiting distribution is obtained as the image of some Poisson random measure.

\subsection{Scale-free random trees}

Scale-free random trees form a one-parameter family of random trees that grow following a preferential attachment algorithm; see \cite{BA}. 
Fix a parameter
$\beta\in(-1, \infty)$, and start  for $n=1$ from  
the unique tree $T^{(\beta)}_1$ on ${\mathcal V}_1=\{0,1\}$ which has a single edge connecting $0$ and $1$. Then suppose that $T^{(\beta)}_n$ has been constructed for some $n\geq 1$, and for every $i\in{\mathcal V}_n=\{0,\ldots,n\}$, denote by $d_n(i)$ the degree of the vertex $i$ in $T^{(\beta)}_n$. Conditionally given $T^{(\beta)}_n$, we construct the tree $T^{(\beta)}_{n+1}$  by incorporating  the new vertex $n+1$  to $T^{(\beta)}_n$ and adding an edge  between $n+1$ and a vertex $v_n\in\{0,\ldots,n\}$ chosen at random according to the law 
$$\P(v_n=i)= \frac{d_n(i)+\beta}{2n+\beta(n+1)}\,, \qquad i\in\{0,\ldots,n\}.$$
Recall that there is the identity $\sum_{i=0}^nd_n(i)=2n$ (because $T^{(\beta)}_n$ is a tree with $n+1$ vertices), so  the preceding indeed defines a probability on $\{0,\ldots,n\}$.  Note also that when one let $\beta\to \infty$,  then $v_n$ becomes uniformly distributed on $\{0,\ldots,n\}$, and  the algorithm yields a uniform recursive tree as in the preceding section.

Just as for recursive trees, one can check that the conditions \eqref{Hk} are fulfilled for all $k\geq 1$ with $\ell(n)=\ln n$ and $L_k\equiv k(1+\beta)/(2+\beta)$; see for instant Section 4.4 in \cite{Durrett} in the case $\beta=0$. Hence we know from Theorem \ref{T1} and Proposition \ref{P1}  that percolation  in the regime \eqref{E1} produces a single giant cluster, and more precisely that \eqref{E2} holds with $G(c) \equiv \theta(c) = \e^{- c(1+\beta)/(2+\beta)}$. 
It has been shown recently in \cite{BU} that asymptotic behavior of the sizes of the largest clusters for percolation on a scale-free tree is similar to that on a random recursive tree.

\begin{theorem}\label{T3} Let $T_n=T^{(\beta)}_n$ denote a random scale free tree on $\{0,1, \ldots, n\}$ with parameter $\beta > -1$. For every fixed integer $j$,  in the regime \eqref{E1} with $\ell(n)=\ln n$, 
$$\left(\frac{\ln n}{n}C^1_{p(n)}, \ldots, \frac{\ln n}{n}C^j_{p(n)}\right)$$
converges in distribution towards
$$(\x_1, \ldots, \x_j)$$
where $\x_1>\x_2>\ldots$ denotes the sequence of the atoms of a Poisson random measure
on $(0,\infty)$ with intensity $c\e^{- c(1+\beta)/(2+\beta)} x^{-2}\d x\,.$
\end{theorem}

The key splitting property of random recursive trees fails for scale-free random trees, and the approach in \cite{BU} for establishing Theorem \ref{T3} thus departs significantly from that for Theorem \ref{T2}. In short, one  superposes Bernoulli bond percolation to the growth algorithm with preferential attachment as follows. Each time an edge is inserted, 
we draw an independent Bernoulli variable $\epsilon$ with parameter $p(n)$. If $\epsilon=1$, the edge is left intact, otherwise we  cut this edge in two at its mid-point.
The upshot of cutting  rather than removing edges is that the former procedure preserves the degrees of vertices, where the degree of a vertex is defined as the sum of the intact edges and half-edges attached to it. This is crucial for running the construction with preferential attachment. 

This  enables us to adapt a classical  idea in this area (see, e.g. \cite{Durrett}), namely to consider a continuous time version of the growth algorithm with preferential attachment and interpret the latter in terms of a continuous time branching processes. Roughly speaking, incorporating percolation to the algorithm yields systems of branching processes with rare neutral mutations, where a mutation event corresponds to the insertion of an edge that is cut in its mid-point. Each branching process in the system corresponds to a percolation cluster which grows following a dynamic with preferential attachment.  One has to study carefully the asymptotic behavior of such systems of branching processes with neutral mutations, and then derive Theorem \ref{T3}.

\end{section}

\noindent{\bf Acknowledgments}.  I would like to thank  David Croydon for a question that he raised during the workshop {\it Random Media II} at Tohoku University in September 2012, which has  motivated Section 2  of the present text

  \end{document}